\documentclass{amsart}
\usepackage{amsmath}%
\usepackage{amsfonts}%
\usepackage{amssymb}%
\usepackage{graphicx}
\theoremstyle{plain}
\newtheorem{acknowledgement}{Acknowledgement}

\newtheorem{conclusion}{Conclusion}

\newtheorem{definition}{Definition}

\numberwithin{equation}{section}
\begin{document}
\title[Distances of MA process error terms]{Asymptotic probability distribution of distances between local extrema of error terms of a moving average process}
\author{Argyn Kuketayev}
\email[A. Kuketayev]{jawabean@gwmail.gwu.edu}%
\date{May 22, 2011}

\begin{abstract}
Consider error terms $\xi_i$ of a moving average process MA(q), where $\xi_i=\sum_{j=0}^q{\varepsilon_{i-j}}$ and $\varepsilon_i$ - independent identically distributed (i.i.d.) random variables.
We recognize a term $\xi_i$ as a local maximum if the following condition holds true: $ \xi_{i-1} < \xi_i > \xi_{i+1}$.
If the local maximum $\xi_i$ is followed by the next local maxiumum $\xi_k$, then $d=k-i$ is the distance between local maxima. The distances $d_j$ themselves are random vriables. In this paper we study the probability distribution of distances $d_j$. Particularly, we show that for any $q>0$ mean distance $E[d_j]=4$ and asymptotically the variance is also equal to 4.
\end{abstract}
\maketitle

\section{Average distance between local maxima}

\begin{definition}
Error terms of a certain MA(q) process are given by equation $\xi_i=\sum_{j=0}^q{\varepsilon_{i-j}}$, where $\varepsilon_i$ - i.i.d. random variables, see \cite{Hamilton}.
\end{definition}

\begin{definition}
A term $\xi_i$ is a local maximum (peak), if the following condition is true $ \xi_{i-1} < \xi_i > \xi_{i+1}$.
\end{definition}

Expanding the definition of a local maximum by susbstituting the MA(q) error term expressions, we can 
re-write the local maxiumum condition as follows: 
\[
\sum_{j=0}^q{\varepsilon_{i-j-1}} < \sum_{j=0}^q{\varepsilon_{i-j}} > \sum_{j=0}^q{\varepsilon_{i-j+1}}
\]
This equation breaks down to independent conditions:
\[
\varepsilon_{i-q-1} < \varepsilon_{i} 
\]
and
\[
\varepsilon_{i-q} > \varepsilon_{i+1}
\]
for any $q>0$.
Hence, the probability to encounter a local maximum at a term $\xi_i$ can be computed as
\[
Pr[max]=Pr[\varepsilon_{i-q-1} < \varepsilon_{i}]\cdot Pr[\varepsilon_{i-q} > \varepsilon_{i+1}]
\]

We shall use a standard cumulative distribution function (CDF) $F(\varepsilon)$, defined as a probability of $\varepsilon\leq \varepsilon_{i}$: \[F(\varepsilon_{i})=Pr[\varepsilon\leq \varepsilon_{i}] = \int^{\varepsilon_i}_{-\infty}f(\varepsilon)\cdot d\varepsilon\]
, where $f(\varepsilon)$ is PDF (probability density function) of $\varepsilon$. This also could be written as 
\[F(\varepsilon_{i}) = \int^{F(\varepsilon_i)}_{0}dF(\varepsilon)\]
On the other hand, probability of $\varepsilon>\varepsilon_{i}$ is \[\int^{1}_{F(\varepsilon_i)}dF(\varepsilon)=1-F(\varepsilon_{i})\]
 
Now, we can find the probabilities
\[ 
Pr[\varepsilon_{i-q-1} < \varepsilon_{i}]
 = \int^{1}_{0}dF(\varepsilon_{i-q-1})\cdot \int^{1}_{F(\varepsilon_{i-q-1})}dF(\varepsilon_{i})
 =\frac{1}{2}
 \]
\[ 
Pr[\varepsilon_{i-q} < \varepsilon_{i+1}]
 = \int^{1}_{0}dF(\varepsilon_{i-q})\cdot \int_{0}^{F(\varepsilon_{i-1})}dF(\varepsilon_{i})
 =\frac{1}{2}
 \]
 then the probability of a local maximum is

\begin{equation}
Pr[max]=\frac{1}{2}\cdot\frac{1}{2}=\frac{1}{4}
\end{equation}
 
Subsequently, in the set of $N$ error terms $\xi_i={\xi_1,\xi_2,...,\xi_N}$ the expected number of 
peaks is equal to $Pr[max]\cdot N$, then the estimate of the mean distance between them is 
\begin{equation}
\tilde{d}= \lim_{N \rightarrow \infty}{\frac{N}{N\cdot Pr[max]-1}} = 
 \lim_{N \rightarrow \infty}{\frac{1}{Pr[max]-\frac{1}{N}}} = 
4
\end{equation}
 
\section{Probability distribution of distances between local maxima}

\subsection{Probability of distance equal to 2}

The following condition must hold true in order to recognize after the peak $\xi_i$ the nearest peak at $\xi_{i+2}$, i.e. at the distance equal to 2:
\[
\xi_{i-1} < \xi_i > \xi_{i+1} < \xi_{i+2} > \xi_{i+3}
\]
This equation is similar to the conditions defining the local maxima in the underlying sequence $\varepsilon_i$. The distances between maxima of $\varepsilon_i$ were studied in \cite{argyn}, and a similar problem was studied in \cite{oshanin} while generating the random sequences from permutations. Unlike the results of the mentioned two works, MA(q) process studied in this paper yields much simpler equations for the probabilities of distances between maxima. It is caused by \itshape{separation} \normalfont  of multiple integrals, which define the probabilities, into products of trivial 2nd-order integrals.

Let us expand the condition using the MA(q) process' error terms defintion:
\[
\sum_{j=0}^q{\varepsilon_{i-j-1}} 
< \sum_{j=0}^q{\varepsilon_{i-j}} 
> \sum_{j=0}^q{\varepsilon_{i-j+1}}
< \sum_{j=0}^q{\varepsilon_{i-j+2}}
> \sum_{j=0}^q{\varepsilon_{i-j+3}}
\]

If $q>2$ then this condition breaks down into the following independent conditions:
\[
\varepsilon_{i-q-1} < \varepsilon_{i} 
\]
\[
\varepsilon_{i-q} > \varepsilon_{i+1}
\]
\[
\varepsilon_{i-q+1} > \varepsilon_{i+2}
\]
\[
\varepsilon_{i-q+2} > \varepsilon_{i+3}
\]

The joint probability of satisfying these conditions is equal to:

\[
Pr[\xi_{i-1} < \xi_i > \xi_{i+1} < \xi_{i+2} > \xi_{i+3}] =
\]
\[
=Pr[\xi_{i-1} < \xi_i ] \cdot
Pr[ \xi_i > \xi_{i+1} ] \cdot
Pr[ \xi_{i+1} < \xi_{i+2}] \cdot
Pr[\xi_{i+2} > \xi_{i+3}] 
\]

As we have shown above, all these four probabilities are equal to $\frac{1}{2}$.

In order to compute the probability $Pr[d=2]$ of the distance between local maxima equal to 2, we shall use Bayes equation for conditional probabilities
\begin{equation}
Pr[A|B]=\frac{Pr[A\cap B] }{Pr[B]}
\end{equation}, where an event $A\cap B$
 is sequence of error terms $\xi_i={\xi_1,\xi_2,...,\xi_N}$ with one maximum in its head and one maximum in its tail, event B is the maximum in first three terms, and event $A|B$ is two local maxima on a given distance from each other. 
Hence,
\begin{equation}
Pr[d=2]=\frac{Pr[\xi_{i-1} < \xi_i > \xi_{i+1} < \xi_{i+2} > \xi_{i+3}] }{Pr[max]}
=\frac{\frac{1}{2^4} }{\frac{1}{4}}=\frac{1}{4}
\end{equation}

\subsection{Probability of distance equal to 3}

The following condition must hold true in order to recognize the nearest peak $\xi_{i+3}$ from the peak $\xi_i$:
\[
\left(
\xi_{i-1} < \xi_i > \xi_{i+1} < \xi_{i+2} < \xi_{i+3} > \xi_{i+4} 
\right)
\bigcup
\left(
\xi_{i-1} < \xi_i > \xi_{i+1} > \xi_{i+2} < \xi_{i+3} > \xi_{i+4} 
\right)
\]

Similar to distance $d=2$, this condition can be expanded using the MA(q) process' error terms defintion,
then can be broken down to the following independent conditions when $q>3$:
\[
\varepsilon_{i-q-1} < \varepsilon_{i} 
\]
\[
\varepsilon_{i-q} > \varepsilon_{i+1}
\]
\[
\left(
\varepsilon_{i-q+1} < \varepsilon_{i+2}
\right)
\bigcup
\left(
\varepsilon_{i-q+1} > \varepsilon_{i+2}
\right)
\]
\[
\varepsilon_{i-q+2} < \varepsilon_{i+3}
\]
\[
\varepsilon_{i-q+3} > \varepsilon_{i+4}
\]

As we have shown above, all the four probabilities are equal to $\frac{1}{2}$. The expression on the third line 
can be evaluated as follows:
\[
Pr[
\left(
\varepsilon_{i-q+1} < \varepsilon_{i+2}
\right)
\bigcup
\left(
\varepsilon_{i-q+1} > \varepsilon_{i+2}
\right)
]=\frac{1}{2}+\frac{1}{2}=1
\]

Now, the probability $Pr[d=3]$ of distance $d=3$ is given by equation:
\[
Pr[d=3]=\frac{ \frac{1}{2} \cdot \frac{1}{2} \cdot 1 \cdot \frac{1}{2} \cdot \frac{1}{2} }
{Pr[max]}
= \frac{ \frac{1}{2^4} } { \frac{1}{4} }
= \frac{1}{4}
\]
\subsection{Asymptotic probability distribution of distances}

Provided that $q>d$, it can be easily shown that for any distance $d$ the probability $Pr[d]$ 
can be computed as follows:
\[
Pr[d]=\frac{ \frac{1}{2} \cdot \frac{1}{2} \cdot \pi(d) \cdot \frac{1}{2} \cdot \frac{1}{2} }
{Pr[max]}
= \frac{\pi(d)}{4}
\]
\[
\pi(d)=\frac{ d-1 }{2^{d-2}}
\],
where $\pi(d)$ is the probability to not encounter a local maxima in the sequence
$\xi_{i+2},\xi_{i+3},...,\xi_{i+d-2}$.

Finally, by combining the equations for probabilities of all possible distances $d$, we can
write the probability mass function (PMF) of the distribution of the distances between local maxima of error terms of this MA(q) process as
follows:
\begin{equation}
Pr[d]= \frac{d-1}{2^d}
\end{equation}
where $q>d$.

The asymptotic estimate of a mean distance is
\[
E[d]=\sum_{d=2}^{\infty}{ \lim_{d<q \rightarrow \infty}\frac{d-1}{2^d} d}=4
\] 
and of the variance is
\[
Var[d]=\sum_{d=2}^{\infty}{ \lim_{d<q \rightarrow \infty}{\frac{d-1}{2^d}}d^2}-E[d]^2=4
\] 

The PMF $Pr[d]$ decreases exponentially, therefore even for MA(10) process, the asymptotic PMF is a reasonably good approximation of the exact PMF. The exact PMF for MA(1) was obtained in \cite{zubkov}.

\begin{conclusion}
blah-blah
\end{conclusion}
\begin{acknowledgement}
Author is very gratefull to Dr. F.M.Pen'kov for fruitfull discussions and interesting suggestions on this topic.
\end{acknowledgement}

\end{document}